\documentclass[a4paper]{amsart}
\usepackage[english]{babel}
\usepackage[utf8x]{inputenc}
\usepackage[T1]{fontenc}
\usepackage{amsthm}
\usepackage{dsfont}
\usepackage{amssymb}
\usepackage{amsmath}
\theoremstyle{plain}
\usepackage{chngcntr}
\usepackage{relsize}
\usepackage{pdfpages}

\newtheorem{Lemma}{Lemma}

\newtheorem{Theorem}[Lemma]{Theorem}
\newtheorem{Proposition}[Lemma]{Proposition}
\newtheorem{Corollary}[Lemma]{Corollary}

\newtheorem{remark}{Remark}

\usepackage[a4paper,top=3cm,bottom=2cm,left=3cm,right=3cm,marginparwidth=1.75cm]{geometry}

\usepackage{float}
\usepackage{nccmath}
\usepackage{mathtools}
\usepackage{amsfonts}
\usepackage{amsmath}
\usepackage{graphicx}
\usepackage[colorinlistoftodos]{todonotes}
\usepackage[colorlinks=true, allcolors=blue]{hyperref}
\usepackage{thmtools}
\title{Character sums over co-dimension one sub-lattices in finite field extensions}
\author{Aishik Chattopadhyay}
\address{Aishik Chattopadhyay,
Ramakrishna Mission Vivekananda Educational and Research Institute, Department of Mathematics, G. T. Road, PO Belur Math, Howrah, West Bengal 711202, India}
\email{aishik.ch@gmail.com}
\begin{document}
\begin{abstract}
We obtain nontrivial bounds for multiplicative character sums over codimension-one sublattices of finite field extensions $\mathbb{F}_{p^d}$. This extends earlier results of Davenport--Lewis and Chang from full-dimensional settings to the codimension-one setting. As an application, we establish cancellation in character sums over binary cubic forms for intervals of length $p^{3/8+\varepsilon}$. Our approach combines Burgess-type amplification technique, representation of character sum using multiplicative energy and it's estimates to prove the results.
\end{abstract}
\maketitle
\tableofcontents
\section{Introduction}
Character sums over short intervals play a central role in analytic number theory. By a character sum over a short interval, we mean a sum of the form 
\[
\sum_{x \in I} \chi(x),
\]
where $\chi$ is a multiplicative character modulo $q$ and $I \subset \mathbb{Z}$ is an interval of length $|I| < q^{1/2}$.

Pólya-Vinogradov inequality states that $\sum_{x \in I} \chi(x)\ll q^{1/2}\log{q}$ where $\chi$ is a character modulo $q$ \cite[Theorem 12.5]{TV}. It obtains a nontrivial upper bound whenever $|I|\geq q^{1/2+\varepsilon}$. A fundamental result in this area is due to Burgess \cite{B1}, who established a nontrivial bound for character sums over intervals of length $|I| > q^{1/4+\varepsilon}$ for any $\varepsilon > 0$. This represents a major improvement over the Pólya--Vinogradov inequality in the short interval regime. Under the Generalized Riemann Hypothesis (GRH), however, nontrivial estimates are known for intervals as short as $|I| > q^{\varepsilon}$, which indicates a substantial gap between conditional and unconditional results. In particular, the Generalized Riemann Hypothesis implies Vinogradov's conjecture that the least quadratic non-residue $n(q)$ satisfies $n(q) \ll q^{\varepsilon}$ which is one of the applications of character sums in a short interval.

In higher dimensions, Davenport and Lewis \cite{DL} initiated the study of character sums over finite field extensions. Let $\mathbb{F}_{p^d}$ be equipped with a basis $\{\omega_1,\dots,\omega_d\}$ over $\mathbb{F}_p$, and consider the linear form 
\[
F(x_1,\dots,x_d) = x_1\omega_1 + \cdots + x_d\omega_d.
\]
They proved that
\[
\sum_{x_1 \in I_1,\dots,x_d \in I_d} \chi(F(x_1,\dots,x_d))
= O\bigl(|I_1|\cdots |I_d|\, p^{-\delta(\varepsilon)}\bigr),
\]
if $I_i$ are intervals satisfying $|I_i| = p^{\rho_d+\varepsilon}$ for all $i$, where $\rho_d = \tfrac{1}{2} - \tfrac{1}{2(d+1)}$. Subsequent work showed that stronger bounds are possible under additional structural assumptions on the basis (see \cite{B2,K}).

Chang \cite{Cha} obtained a significant improvement for the powered basis
\[\{1,\omega,\omega^2,\dots,\omega^{d-1}\},\]
where $\omega$ generates $\mathbb{F}_{p^d}$ over $\mathbb{F}_p$. More precisely, for $d\ge2$, she established a non-trivial estimate when
\[
|I_i|=p^{\rho'_d+\varepsilon} \text{ for all } i,
\text{ with }
\rho'_d=
\frac{\sqrt{d^2+2d-7}+3-d}{8}.
\]
Note that $\rho'_d<\rho_d$. In particular, for $d=2$, Chang proved nontrivial bounds for character sums if the intervals are of length at least $p^{1/4+\varepsilon}$. Without imposing special assumptions on the basis, she later showed in \cite{C1} that for arbitrary $d\ge1$ one still obtains nontrivial estimates for intervals of length $p^{2/5+\varepsilon}$. In particular, this improves her result in \cite{Cha} if $d\geq 5$. Later, this was improved by Konyagin \cite{Kon}, who reduced the exponent to $1/4+\varepsilon$.

 \medskip
The main purpose of this paper is to extend Chang's framework \cite{C1} to sums taken over \emph{a sublattice of codimension one} in $\mathbb{F}_{p^d}$. More precisely, instead of summing over all coordinates, we consider sums of the form
\[
\sum_{x_0,\dots,x_{d-2}} \chi\bigl(x_0 + x_1\omega + \cdots + x_{d-2}\omega^{d-2}\bigr),
\]
which correspond to restricting the summation to a hyperplane in $\mathbb{F}_{p^d}$. To the best of our knowledge, such sums have not previously been studied.

Our main result stands as follows.
\begin{Theorem}[Main Theorem]\label{MTG}
Let $\chi$ be a nontrivial multiplicative character of $\mathbb{F}_{p^d}$, and let $\omega$ generate $\mathbb{F}_{p^d}$ over $\mathbb{F}_p$. For any $\varepsilon > 0$, there exists $\delta(\varepsilon) > 0$ such that if $J_0,\dots,J_{d-2}$ are intervals of length $|J_i| \ge p^{\rho''_d+\varepsilon}$, where
\[
\rho''_d =
\frac{-50d+73+5\sqrt{100d^2-196d+49}}{48},
\]
then
\[
\left| \sum_{x_0 \in J_0,\dots,x_{d-2} \in J_{d-2}}
\chi\bigl(x_0 + x_1\omega + \cdots + x_{d-2}\omega^{d-2}\bigr) \right|
\ll p^{-\delta(\varepsilon)} \prod_{i=0}^{d-2} |J_i|,
\]
where the implied constant is independent of $\omega$.
\end{Theorem}
L. Mérai, I. E. Shparlinski, and A. Winterhof \cite{MIW} studied mixed character sums involving polynomial values over the set of vectors in $\mathbb{F}_{p^d}$ having Hamming weight $s$. They referred to such elements as $s$-sparse elements. More precisely, let $\{\omega_1,\omega_2,\ldots,\omega_d\}$ be an ordered basis of $\mathbb{F}_{p^d}$ over $\mathbb{F}_p$ and for $v=a_1\omega_1+a_2\omega_2+\cdots+a_d\omega_d\in\mathbb{F}_{p^d},$ let the Hamming weight of v be defined as the number of
nonzero components among $a_1,a_2,\ldots, a_d$. They obtained a square-root cancellation of the complete character sum over $s-$sparse elements in $\mathbb{F}_{p^d}$. In Theorem~\ref{MTG}, the elements appearing in the character sum have zero $\omega^{d-1}$-coordinate and hence are at most $(d-1)$-sparse with respect to the ordered basis $\{1,\omega,\ldots,\omega^{d-1}\}$. Thus, our result is related to the estimates for character sums over sparse elements obtained by Mérai, Shparlinski, and Winterhof.

Note that after averaging over an arbitrary set in the $d-$th coordinate, the above character sum still gives a non-trivial estimate, with the resulting saving independent of both the size and the sparsity of the set. In \cite{II}  Iyer and Shparlinski  studied character sums over elements of finite field extensions with sparse coordinates which is also related to sparse character sum taken above. 

As an immediate consequence of Theorem \ref{MTG}, we obtain the following result in dimension $d=3$.
\begin{Corollary}\label{T}
Let $\varepsilon > 0$. Then there exists $\delta(\varepsilon) > 0$ such that for all sufficiently large primes $p$, the following holds. If $\omega$ generates $\mathbb{F}_{p^3}$ over $\mathbb{F}_p$ and $I,J \subset \mathbb{F}_p$ are intervals of length at least $p^{3/8+\varepsilon}$, then for every nontrivial multiplicative character $\chi$ of $\mathbb{F}_{p^3}$,
\[
\left| \sum_{x \in I} \sum_{y \in J} \chi(x + \omega y) \right|
\ll p^{-\delta(\varepsilon)} |I||J|.
\]
\end{Corollary}
Let $f$ be an irreducible polynomial over $\mathbb{F}_p$. The de-homogenization $f(x,1)$ of the form $f(x,y)$ has a root $\omega$ in $\mathbb{F}_{p^3}$. If
$
\operatorname{N}_{\mathbb{F}_{p^3}/\mathbb{F}_p}
:\mathbb{F}_{p^3}\rightarrow \mathbb{F}_p
$
denotes the norm map then
\[
\chi(f(x,y))
=
\chi\!\left(
\operatorname{N}_{\mathbb{F}_{p^3}/\mathbb{F}_p}(x+\omega y)
\right)
=
(\chi\circ
\operatorname{N}_{\mathbb{F}_{p^3}/\mathbb{F}_p})(x+\omega y).
\]
Since $\chi\circ\operatorname{N}_{\mathbb{F}_{p^3}/\mathbb{F}_p}$ is a multiplicative character of $\mathbb{F}_{p^3}$, our next corollary below  follows immediately by applying Corollary~\ref{T} to this character.
\begin{Corollary}\label{MC}
Given $\varepsilon>0$, there is $\delta(\varepsilon)>0$ such that the following holds.
Let $p$ be a large prime and $f(x,y)=x^3+ax^2y+bxy^2+cy^3$ be an irreducible cubic form over $\mathbb{F}_p$. Let $I,J$ be intervals of size at least $p^{3/8+\varepsilon}$. Then
$$\left|\sum_{x\in I,y\in J}\chi(f(x,y))\right|\ll p^{-\delta(\varepsilon)}|I||J|$$
for any nontrivial multiplicative character $\chi(\bmod{p)}$. The implied constant and the saving is uniform in $f$.
\end{Corollary}

Moreover, we can get the same exponent for the size of intervals to obtain a non-trivial character sum on homogeneous binary non-singular cubic forms. Recall that a homogeneous binary form is said to be non-singular if its de-homogenization has no multiple roots. We shall come to this at the end as an remark.

Pierce and Xu \cite{Pichu} established Burgess-type bounds for character sums in arbitrary dimension evaluated at \emph{admissible} homogeneous forms of arbitrary degree via recently developed stratification results by Xu. Their result is the following:
\begin{Theorem}
Let $\chi$ be a non-principal multiplicative character of order $\Delta$
modulo a prime $p$. Let $n\geq 2$ be fixed. For each $r\geq 1$, define
\[
\Psi=\Psi_{n,r}
=\left\lfloor\frac{r-1}{n-1}\right\rfloor.
\]
Let $\mathbf{H}=(H_1,\ldots,H_n)\in\mathbb{R}_{\geq 1}^n$ have maximum
element $H_{\max}$ and minimum element $H_{\min}$, and assume that
$H_{\max}H_{\min}<p^{1+1/(2\Psi)}.$ Then, for every form $F\in\mathbb{Z}[x_1,\ldots,x_n]$ of degree $d$ whose reduction modulo $p$ is \emph{admissible}, uniformly in $\mathbf{N}=(N_1,\ldots,N_n)$ and for every integer $r\geq 1$, we have
\[
\sum_{\substack{\mathbf{x}\in\mathbb{Z}^n\\
x_i\in(N_i,N_i+H_i]}}
\chi(F(\mathbf{x}))
\ll
(H_1\cdots H_n)^{1-\frac{1}{2r}}
H_{\min}^{-\frac{1}{2r}}
p^{(n\Psi+1)/4r\Psi}
(\log p)^{n+1}.
\]
Here the implied constant depends on $D,\Delta,n,r$ and is independent
of $F$ and $\mathbf{N}$.
\end{Theorem}
In particular, for binary admissible forms, the above theorem yields a non-trivial estimate for intervals of length $p^{1/3+\varepsilon}$ (see \cite[equation (1-6)]{Pichu}). In contrast, Corollary~\ref{MC} yields a nontrivial bound for cubic binary forms when the interval lengths are at least $p^{3/8+\varepsilon}$. The proof follows the standard Burgess-type arguments representing the multiplicative energy through character sums, thereby relating the original character sum to the multiplicative energy of a set of tuples, and then estimate this energy. Although our method does not yield a bound as strong as Pierce and Xu \cite{Pichu}, it obtains a non-trivial result by a
modification of Chang's technique, relying on less advanced tools than algebraic geometric tools used in \cite{Pichu}.

\section{Notation and Preliminaries}\label{1stsec}
Let $A$ be a commutative ring with unity. In our article 
$A=\mathbb{F}_{p^d}$ or $\mathbb{F}_p\times \mathbb{F}_p\times \mathbb{F}_p,
$ or $\mathbb{F}_p\times \mathbb{F}_{p^2}$ according to the context. We denote the group of invertible elements of $A$ by $A^*$. Let $\chi$ be a multiplicative character of $A$.

Throughout the sequel, $c>0$ denotes a sufficiently small constant which may vary from line to line.  $\mathbb Z_{>0}$ is the set of positive integers. For $n\in \mathbb N$, let $d(n)$ denote the divisor function. We repeatedly use the classical estimate
\[
d(n)= \exp\left(O\left(\frac{\log n}{\log\log n}\right)\right).
\]

We now fix some notation.

\begin{enumerate}
\item For subsets $B_1,\ldots,B_n\subset A$, define the multiplicative energy
$
E(B_1,\ldots,B_n)
$
to be the number of tuples
$(b_1,\ldots,b_n,b_1',\ldots,b_n')$
such that
\[
b_1\cdots b_n=b_1'\cdots b_n',
\qquad
b_j,b_j'\in B_j\cap A^*.
\]
For $A=\mathbb{F}_{p^d}$,
\[
E(B_1,\ldots,B_n)=\]
\[\left|
\left\{
(b_1,\ldots,b_n,b_1',\ldots,b_n'):
b_1\cdots b_n=b_1'\cdots b_n'\neq 0,\;b_i,b_i'\in B_i \text{ for all } 1\leq i\leq n
\right\}
\right|.
\]
Using the orthogonality relation of multiplicative characters, we obtain
\begin{equation}\label{mcc}
E(B_1,\ldots,B_n)
=
\frac{1}{|A^*|}
\sum_\chi
\prod_{i=1}^n
\left|
\sum_{\eta\in B_i}\chi(\eta)
\right|^2.
\end{equation}
\item \label{sec} We briefly recall the Burgess-type amplification argument that will be used in the paper. Consider the character sum
$
\left|
\sum_{x_1\in J_1, x_2\in J_2,\ldots,x_d\in J_d}
\chi(x_1+\omega x_2+\omega^2 x_3+\ldots +\omega^{d-1}x_d)
\right|,
$
where $\omega$ is a root of an irreducible polynomial of degree d over $\mathbb F_p$, and $J_1,J_2,\ldots,J_d$ are intervals.
Choose auxiliary intervals $J_1^0,J_2^0,\ldots, J_d^0$ together with an integer $l$ satisfying
$l|J_i^0|<p^{-\varepsilon}|J_i|,$ for all $1\leq i\leq d$
for some sufficiently small $\varepsilon>0$. Translating $(x_1,x_2,\ldots, x_d)$ by $(tx_1^0,tx_2^0,\ldots,tx_d^0)$ with
$
t\in [1,l],
\,
(x_1^0,x_2^0,\ldots,x_d^0)\in J_1^0\times J_2^0\times \cdots\times J_d^0,
$ it suffices to estimate
\begin{equation}\label{Bamp}
\frac{1}{l\,|J_1^0|\,|J_2^0|\,\ldots|J_d^0|}
\left|
\sum_{\substack{
1\le t\le l\\
(x_1^0,x_2^0,\ldots,x_d^0)\in I_1^0\times I_2^0\times \cdots I_d^0
}}
\;
\sum_{\substack{x_1\in J_1\\x_2\in J_2\\ \ldots\\x_d\in J_d}}
\chi\bigl(x_1+t x_1^0+(x_2+t x_2^0)\omega+\cdots+(x_d+tx_d^0)\omega^{d-1}\bigr)
\right|.
\end{equation}
The role of the auxiliary intervals $J_1^0, J_2^0,\cdots,J_d^0$ is to create enough averaging to obtain a necessary cancellation in the character sum.
We define
\[
\Phi(\mu)
=
\left|
\left\{
(x_1,x_2,\ldots, x_d,x_1^0,x_2^0,\ldots, x_d^0)\in
\prod_{i=1}^dJ_i\times \prod_{i=1}^d J^0_i:
\mu=\frac{x_1+x_2\omega+\cdots+x_d\omega^{d-1}}{x_1^0+x_2^0\omega+\cdots+x_d^0\omega^{d-1}}
\right\}
\right|.
\]
Therefore, the above quantity \eqref{Bamp} equals
\begin{equation}\label{beq}
\frac{1}{l\,|J_1^0|\,|J_2^0|\cdots|J_d^0| }\sum_{\mu\in\mathbb F_{p^d}}
\Phi(\mu)
\left|
\sum_{1\le t\le l}\chi(t+\mu)
\right|.
\end{equation}
Next, we apply H\"older's inequality in the equation \ref{beq} to obtain
\begin{equation}\label{HE}
\eqref{beq}
\leq\frac{1}{l\,|J_1^0|\,|J_2^0|,\ldots, |J^0_d|}
\left[
\sum_{\mu\in\mathbb F_{p^d}}
\Phi(\mu)
\right]^{1-\frac1k}
\left[
\sum_{\mu\in\mathbb F_{p^d}}
\Phi(\mu)^2
\right]^{\frac1{2k}}
\left[
\sum_{\mu\in\mathbb F_{p^d}}
\left|
\sum_{1\le t\le l}\chi(\mu+t)
\right|^{2k}
\right]^{\frac1{2k}}.
\end{equation}
The first, second, and third factors above can be estimated using the definition of $\Phi(\mu)$, bounds for multiplicative energy, and Weil's bound for character sums, respectively. We shall revisit the details in the proof of our Theorem \ref{MTG} where we shall use this amplification method.
\end{enumerate}
We recall the Weil bound for multiplicative character sums in finite fields, which plays a key role in our estimates.
\begin{Proposition}[Weil bound for multiplicative characters on $\mathbb{F}_{p^d}$]\label{WEL}
Let $q=p^d$ be a power of a prime and let $\chi$ be a nontrivial multiplicative character of $\mathbb{F}_{q}$ of order $n$. Let $f(x)\in\mathbb{F}_q[x]$ be a polynomial which has $m\ge 1$ distinct roots and is not an $n$-th power in $\mathbb{F}_q[x]$. Then 
\[
\left|\sum_{x\in\mathbb{F}_q}\chi\big(f(x)\big)\right|\le (m-1)q^{1/2}.
\]
\end{Proposition}
\begin{proof}
See \cite[Theorem 2C', p.~43]{Sch}.
\end{proof}
The following lemma provides a key estimate used in the main argument.
\begin{Lemma}\label{WLM}
Let $\chi$ be a multiplicative character of $\mathbb{F}_{p^d}$ of order $n>1$, and let $I \subset \mathbb{Z}$ be an interval. Then for every integer $r \ge 1$, we have
\[
\sum_{u \in \mathbb{F}_{p^d}}
\left|\sum_{z \in I} \chi(u+z)\right|^{2r}
\ll
2r\, p^{d/2} |I|^{2r}
\;+\;
r^{2r}|I|^r\,p^d\,,
\]
with an absolute implied constant.
\end{Lemma}

\begin{proof}
Expanding the $2r$-th moment gives
\[
\sum_{u \in \mathbb{F}_{p^d}}
\left|\sum_{z \in I} \chi(u+z)\right|^{2r}
=
\sum_{z_1,\ldots,z_{2r}\in I}
\sum_{u \in \mathbb{F}_{p^d}}
\chi\!\left(
\prod_{i=1}^{r}(u+z_i)
\prod_{i=r+1}^{2r}(u+z_i)^{-1}
\right).
\]
\[
=\sum_{z_1,\ldots,z_{2r}\in I}
\sum_{u \in \mathbb{F}_{p^d}}
\chi\!\left(
\prod_{i=1}^{r}(u+z_i)
\prod_{i=r+1}^{2r}(u+z_i)^{p^d-2}
\right).
\]
The last equality follows as $u^{p^d-2}$ is the inverse of an element $u\in \mathbb{F}_{p^d}^*$. Fix a tuple $\mathbf{z}=(z_1,\ldots,z_{2r}) \in I^{2r}$ and define the polynomial function
\[
F_{\mathbf{z}}(u)
=
\prod_{i=1}^{r}(u+z_i)
\prod_{i=r+1}^{2r}(u+z_i)^{p^d-2}.
\]

We classify $\mathbf{z}\in I^{2r}$ into two types, i.e., bad tuple and good tuple. $\mathbf{z}$ is called \emph{bad tuple} if the polynomial $F_{\mathbf{z}}(u)$ becomes a perfect $n$-th power in $\mathbb{F}_{p^d}[u]$. This can only happen if every distinct value among $z_1,\ldots,z_{2r}$ occurs at least twice in this collection and the tuple contains at most $r$ distinct elements. Hence, the number of bad tuples is at most $ r^{2r}|I|^r$. For each bad tuple, we use the trivial bound $p^d$ for the character sum. So, their total contribution is $r^{2r}|I|^rp^d.$ For all remaining tuples (which are called \emph{good tuples}), the polynomial $F_{\mathbf{z}}(u)$ is not a perfect $n$-th power in $\mathbb{F}_{p^d}[u]$. Hence, by Proposition \ref{WEL}, the character sum due to these tuples is estimated by $2r\, p^{d/2}$. Since there are at most $|I|^{2r}$ tuples in total, the contribution of good tuples is
$
\ll 2r\, p^{d/2} |I|^{2r}.
$
Combining the contributions of bad and good tuples completes the proof.
\end{proof}
\section{Multiplicative energy of intervals in $\mathbb{F}_{p^d}$}
In the previous section, we mentioned that estimating the character sum requires bounding certain multiplicative energies. We begin with the following lemma, which will be used in the energy estimates.
\begin{Lemma}\label{GL}
Let $\omega \in \mathbb{F}_{p^d}$ be a generator of $\mathbb{F}_{p^d}$ on $\mathbb{F}_p$. Define
\[
Q=\left\{x_0+x_1\omega+\cdots+x_{d-2}\omega^{d-2}
:\;
x_i\in[1,a]
\text{ for all }0\le i\le d-2
\right\},
\]
and
\[
Q_0=
\left\{
y_0+y_1\omega
:\;
y_i\in[1,b]
\text{ for } i=0,1
\right\},
\]
where $ab<p/2$. Then, for every $\eta\in\mathbb{F}_{p^d}$,
\[
\bigl|
\{(z_1,z_2)\in Q\times Q_0:\eta=z_1z_2\}
\bigr|
=
\exp\!\left(
O\left(\frac{\log p}{\log\log p}\right)
\right).
\]
\end{Lemma}

\begin{proof}
Fix $\eta \in \mathbb{F}_{p^d}$ and suppose that
\[
z_1 z_2 = z_1' z_2' = \eta,
\qquad
z_1,z_1' \in Q,
\quad
z_2,z_2' \in Q_0.
\]
We write
$
z_1 = \sum_{i=0}^{d-2} x_i \omega^i,
z_1' = \sum_{i=0}^{d-2} x_i' \omega^i,$
and
$
z_2 = y_0 + y_1 \omega,
z_2' = y_0' + y_1' \omega.
$
Expanding the products and comparing coefficients with respect to the basis
$\{1,\omega,\dots,\omega^{d-1}\}$, we obtain
\[
\begin{cases}
y_1 x_{d-2} = y_1' x_{d-2}', \\
y_0 x_{d-2} + y_1 x_{d-3}
=
y_0' x_{d-2}' + y_1' x_{d-3}', \\
\hspace{2cm} \vdots \\
y_0 x_1+y_1x_0
=
y_0' x_1'+y_1'x_0', \\
y_0 x_0 = y_0' x_0'.
\end{cases}
\]
Since all variables lie in intervals of positive integers and every monomial above is bounded by $ab<p/2$, each equality holds in $\mathbb Z$.

The first and last equations imply
$
y_1 x_{d-2} = y_1' x_{d-2}'$, $
y_0 x_0 = y_0' x_0'.
$
Hence, by the divisor bound, the number of possibilities for
$
(y_0',y_1',x_0',x_{d-2}')
$
is at most
\[
\exp\!\left(
O\left(
\frac{\log p}{\log\log p}
\right)
\right).
\]
Once these variables are fixed, the remaining coordinates
$
x_1',x_2',\dots,x_{d-3}'
$
are uniquely determined recursively from the above system of equations. 
\end{proof}
As a consequence, we obtain an estimate of multiplicative energy of two \emph{intervals} in $\mathbb{F}_{p^d}$.
\begin{Corollary}\label{GLC}
Let $Q,Q_0$ be as in Lemma~\ref{GL}. Then
\[
   E(Q,\,Q_0)\;= \;
   \exp\!\left(O\left(\frac{\log p}{\log\log p}\right)\right)\,|Q|\,|Q_0| \;,
\]
\end{Corollary}
\begin{proof}
The method follows from the estimate
\[
E(Q,Q_0)\leq |Q||Q_0|\max_{\eta\in \mathbb{F}_{p^d}}\bigl|
\{(z_1,z_2)\in Q\times Q_0:\eta=z_1z_2\}
\bigr|
\]
and the Lemma above.
\end{proof}
Now we require a lemma where we consider the multiplicative energy of a set tuple and it will be used to prove our Theorem \ref{MTG} in the next section.
\begin{Lemma}\label{nL}
Let $0<\sigma<1$, and define
\[
Q=\left\{x_0+x_1\omega+\cdots+x_{d-2}\omega^{d-2}:x_i\in[1,p^\sigma]\right\},
\qquad
Q_0=\left\{y_0+y_1\omega:y_i\in\left[1,\frac{p^{\frac{9}{20}-\frac{4}{20}\sigma}}{2}\right]\right\}.
\]
Further, let $I_s=[1,p^{1/k_s}], \text{ for } s=1,\ldots,r,$
where $k_s\in \mathbb Z_{>0}$ and
$
\sum_{s=1}^r\frac1{k_s}<1.
$
Then for some $c>0$ we have
\[
E(I_1,\ldots,I_r,Q,Q_0)
\ll
\exp\left(c\frac{\log p}{\log\log p}\right)
|Q||Q_0|
\prod_{s=1}^r |I_s|^{\frac1{10}(19-4\sigma)}.
\]
\end{Lemma}

\begin{proof}
By the character-sum representation of multiplicative energy,
\[
E(I_1,\ldots,I_r,Q,Q_0)
=
\frac1{p^d-1}
\sum_\chi
\prod_{s=1}^r
\left|
\sum_{t\in I_s}\chi(t)
\right|^2
\left|
\sum_{\eta\in Q}\chi(\eta)
\right|^2
\left|
\sum_{\eta\in Q_0}\chi(\eta)
\right|^2.
\]
as the sets $I_1,\ldots,I_r,Q,Q_0$ lie in $\mathbb{F}_{p^d}$.
Applying H\"older's inequality, we obtain
\begin{equation}\label{eq:h1}
\begin{aligned}
E(I_1,\ldots,I_r,Q,Q_0)
\leq
&\prod_{s=1}^r
\left(\underbrace{\frac{1}{p^d-1}
\sum_\chi
\left|
\sum_{t\in I_s}\chi(t)
\right|^{2k_s}
\left|
\sum_{\eta\in Q}\chi(\eta)
\right|^2
\left|
\sum_{\eta\in Q_0}\chi(\eta)
\right|^2}_{A_s}
\right)^{1/k_s}
\\
&\times
\left(\underbrace{
\frac{1}{p^d-1}
\sum_\chi
\left|
\sum_{\eta\in Q}\chi(\eta)
\right|^2
\left|
\sum_{\eta\in Q_0}\chi(\eta)
\right|^2}_{B}
\right)^{1-\sum_s1/k_s}.
\end{aligned}
\end{equation}
 Next, we estimate $A_s$ and $B$.
By Corollary~\ref{GLC} with $a=p^{\sigma}, \, b=p^{(\frac{9}{20}-\frac{4}{20}\sigma)}/2$,
\begin{equation}\label{eq:B}
B
=
E(Q,Q_0)
\ll
\exp\left(c\frac{\log p}{\log\log p}\right)
|Q||Q_0|.
\end{equation}
Taking trivial bound over $Q_0$ and using the interpretation of multiplicative energy via character sums,
\begin{equation}\label{eq:As}
\begin{aligned}
A_s
&\le
|Q_0|^2
E(\underbrace{I_s,\ldots,I_s}_{k_s\text{-times}},Q)
\ll
|Q_0|^2
\exp\left(c_{k_s}\frac{\log p}{\log\log p}\right)
E(\mathbb F_p,Q),
\end{aligned}
\end{equation}
where the factor involving $c_{k_s}$ follows from bounding the elements of $\mathbb{F}_p^*$ into $k_s$ factors by divisor bound. We now estimate $E(\mathbb F_p,Q)$. We write
\[
Q=\{x_0+x_1\omega+x_2\omega^2+\cdots+x_{d-2}\omega^{d-2}:x_0,x_1,\ldots,x_{d-2}\in K\},
\]
where
$K=[1,p^\sigma].$ Then
\[
\begin{aligned}
E(\mathbb F_p,Q)
&=
\bigl|
\{(t_1,t_2,\zeta_1,\zeta_2)\in \mathbb F_p^2\times Q^2:
t_1\zeta_1=t_2\zeta_2\neq0
\}
\bigr|
\\
&=
\Bigl|
\Bigl\{
(t_1,t_2,\mathbf x,\mathbf x')
\in
\mathbb F_p^2\times K^{2(d-1)}
:
t_1\sum_{i=0}^{d-2}x_i\omega^i
=
t_2\sum_{i=0}^{d-2}x_i'\omega^i
\neq0
\Bigr\}
\Bigr|.
\end{aligned}
\]
Comparing coefficients of above gives the system
\begin{equation}\label{eq:coeff}
t_1x_i=t_2x_i',
\qquad
i=0,1,\ldots,d-2,
\end{equation}
with at least one non-zero relation. Fix an index $i$ such that $t_1x_i=t_2x_i'\neq0$. Without loss of generality, assume $i=0$. Then \eqref{eq:coeff} implies
\[
t_1=t_2\frac{x_0'}{x_0},
\qquad
\frac{x_j}{x_0}=\frac{x_j'}{x_0'},
\quad
j=1,\ldots,d-2.
\]
In particular, for $j=1$ the number of quadruples
$
(x_0,x_0',x_1,x_1')
$
satisfying
$
x_0x_1'=x_1x_0'
$
is bounded by
\[
E([1,p^\sigma],[1,p^\sigma])
\ll p^{2\sigma}\log p
\]
by \cite{JH}. Once these variables are fixed, there are at most $p$ choices for $t_2$, while each of
$
x_2,\ldots,x_{d-2}
$
admits at most $p^\sigma$ possibilities. Consequently,
\[
E(\mathbb F_p,Q)
\ll
(d-2)\,p^{1+\sigma(d-1)}\log p.
\]
where $d-2$ comes as $i$ lies in between $0$ and $(d-2)$ in equation \eqref{eq:coeff}. Absorbing logarithmic factors into the exponential term, \eqref{eq:As} becomes
\begin{equation}\label{eq:As2}
A_s
\ll
|Q_0|^2
\exp\left(c_{k_s}\frac{\log p}{\log\log p}\right)
p^{1+\sigma(d-1)}.
\end{equation}
Combining \eqref{eq:h1}, \eqref{eq:B}, and \eqref{eq:As2} it yields
\[
\begin{aligned}
E(I_1,\ldots,I_r,Q,Q_0)
\ll\ &
\exp\left(c\frac{\log p}{\log\log p}\right)
|Q||Q_0|
\\
&\times
|Q_0|^{\sum_s1/k_s}
|Q|^{-\sum_s1/k_s}
p^{(1+\sigma(d-1))\sum_s1/k_s}.
\end{aligned}
\]
Since
$
|Q|=p^{(d-1)\sigma}$ and $
|Q_0|=p^{\frac1{10}(9-4\sigma)},
$
we obtain
$
|Q_0|\,|Q|^{-1}p^{1+\sigma(d-1)}
=
p^{\frac1{10}(19-4\sigma)}.
$
Hence
\begin{align*}
E(I_1,\ldots,I_r,Q,Q_0)
\ll&
\exp\left(c\frac{\log p}{\log\log p}\right)
|Q||Q_0|
p^{\frac1{10}(19-4\sigma)\sum_s1/k_s}\\
=&\exp\left(c\frac{\log p}{\log\log p}\right)
|Q||Q_0|\prod_{s=1}^r
|I_s|^{\frac1{10}(19-4\sigma)}
\end{align*}
since $|I_s|=p^{1/k_s}$.
It completes the proof of the Lemma.
\end{proof}

\section{Proof of Theorem \ref{MTG}}
We now prove Theorem \ref{MTG}. The proof uses a standard application of Burgess amplification together with multiplicative energy estimates as in \cite{Cha}. However, the choice of sizes of auxiliary intervals at the beginning of the proof is very delicate to fit our purpose of Burgess amplification and get the resulting size of the intervals, which is optimal.

\begin{proof}
We assume throughout that
$\rho\geq \frac{3}{8}+2\varepsilon.$
Define
$Q_0=\left\{y_0+y_1\omega:y_i\in\left[1,\frac{p^{\frac{9}{20}-\frac{4\rho}{20}}}{2}\right]\right\}.$
Choose integers \(k_1,\dots,k_r\in \mathbb N\) such that
\begin{equation}\label{choi}
\frac65\rho-\frac{9}{20}-2\varepsilon
<\sum_{s=1}^r\frac1{k_s}
<
\frac65\rho-\frac{9}{20}-\varepsilon,
\end{equation}
where \(\varepsilon>0\) is sufficiently small and $r$ depends upon $\varepsilon$. Let
$I=[1,p^{\varepsilon/2}]$ and $
I_s=[1,p^{1/k_s}]
\text{ for all }1\le s\le r.
$ The sizes of the auxiliary intervals \(Q_0\), \(I\), and \(I_s\), together with the choice in \eqref{choi}, are selected so that the amplification of \(Q\) by \(I\cdot\prod_{s=1}^r I_s\cdot Q_0\) is admissible; more precisely,
\[
|I|\prod_{s=1}^r |I_s|p^{9/20-4\rho/20}
<p^{-\varepsilon}|J_i|,
\qquad 0\leq i<2.
\]
These choices also minimize the admissible value of \(\rho\), which will yield the exponent \(\rho_d''\) at the end of the proof. We now consider the character sum
\[
\sum_{x\in Q}\chi(x),
\]
where
\[
Q=
\left\{
x_0+x_1\omega+\cdots+x_{d-2}\omega^{d-2}
:
x_i\in J_i
\right\}.
\]
Translating \(Q\) by
$
I\cdot \prod_{s=1}^r I_s\cdot Q_0
$
we proceed with the Burgess amplification method as described in section \ref{1stsec}. Therefore, it suffices to estimate
\[
\frac1{|I|\prod_{s=1}^r|I_s||Q_0|}
\sum_{\substack{
q\in Q_0\\
(s_1,\dots,s_r)\in\prod_s I_s\\
(x_0,\dots,x_{d-2})\in\prod_i J_i
}}
\left|
\sum_{t\in I}
\chi\!\left(
t+
\frac{x_0+x_1\omega+\cdots+x_{d-2}\omega^{d-2}}
{s_1\cdots s_r\,q}
\right)
\right|.
\]

Define
\[
\Phi(\mu)
=
\left|
\left\{
(x_0,\dots,x_{d-2},s_1,\dots,s_r,q):
\mu=
\frac{x_0+x_1\omega+\cdots+x_{d-2}\omega^{d-2}}
{s_1\cdots s_r\,q}
\right\}
\right|.
\]
Then the above expression equals
\[
\frac1{|I|\prod_{s=1}^r|I_s||Q_0|}
\sum_{\mu\in\mathbb F_{p^d}}
\Phi(\mu)
\left|
\sum_{t\in I}\chi(\mu+t)
\right|.
\]
Applying H\"older's inequality, we obtain
\[
\begin{aligned}
&
\frac1{|I|\prod_{s=1}^r|I_s||Q_0|}
\sum_{\mu\in\mathbb F_{p^d}}
\Phi(\mu)
\left|
\sum_{t\in I}\chi(\mu+t)
\right|
\\
\le&
\frac1{|I|\prod_{s=1}^r|I_s||Q_0|}
\underbrace{\left(
\sum_{\mu\in \mathbb{F}_{p^d}}\Phi(\mu)
\right)^{1-\frac1k}}_{\alpha}\underbrace{\left(
\sum_{\mu\in \mathbb{F}_{p^d}}\Phi(\mu)^2
\right)^{\frac1{2k}}}_{\beta}\underbrace{\left(
\sum_{\mu\in \mathbb{F}_{p^d}}
\left|
\sum_{t\in I}\chi(\mu+t)
\right|^{2k}
\right)^{\frac1{2k}}}_{\gamma}.
\end{aligned}
\]
Varying over all possible elements of $Q_0,I_s$ and $J_i$ in the definition of $\Phi(\mu)$ we get,
\begin{equation}\label{a1}
\alpha=
\left(
|Q||Q_0|p^{\sum_s1/k_s}
\right)^{1-\frac1k}.
\end{equation}
Following the definition of multiplicative energy together with Lemma~\ref{nL} ($\sigma=\rho$) we deduce,
\begin{equation}\label{b1}
\begin{aligned}
\beta
&=
E(Q,Q_0,I_1,\dots,I_r)^{\frac1{2k}}\ll
\left(
\exp\!\left(
c\frac{\log p}{\log\log p}
\right)
|Q||Q_0|
\prod_{s=1}^r
|I_s|^{\frac1{10}(19-4\rho)}
\right)^{\frac1{2k}}.
\end{aligned}
\end{equation}
Further, Lemma~\ref{WLM} yields
\begin{equation}\label{c1}
\gamma
\ll
p^{\frac{\varepsilon}{4}+\frac{d}{2k}}
+
p^{\frac{\varepsilon}{2}+\frac{d}{4k}}.
\end{equation}
It follows from \eqref{a1}, \eqref{b1}, and \eqref{c1} that
\[
\begin{aligned}
\sum_{q\in Q}\chi(q)
&\ll
\frac1{|I|\prod_s|I_s||Q_0|}
\left(
|Q||Q_0|p^{\sum_s1/k_s}
\right)^{1-\frac1k}
\\
&\qquad\times
\left(|Q||Q_0|p^{\frac{(19-4\rho)}{10}\sum_{s}\frac{1}{k_s}}\right)^{\frac1{2k}}
\left(
p^{\frac{\varepsilon}{4}+\frac{d}{2k}}
+
p^{\frac{\varepsilon}{2}+\frac{d}{4k}}
\right).
\end{aligned}
\]
Since
$
|I|=p^{\varepsilon/2},
$ $
\prod_s|I_s|=p^{\sum_s1/k_s},
$ and $|Q_0|\asymp p^{\frac{9}{10}-\frac{4\rho}{10}}$
we obtain
\[
\sum_{q\in Q}\chi(q)
\ll|Q|^{1-1/2k}\left(p^{(4/20k)\rho-(9/20k)-\left((1+4\rho)/20k\right)\sum_{s}\frac{1}{k_s}}\right)(p^{-\varepsilon/4+d/2k}+p^{d/4k}).
\]
Choosing \(k>d/\varepsilon\) we get
\[
\sum_{q\in Q}\chi(q)
\ll |Q|^{1-1/2k}\left(p^{(4/20k)\rho-(9/20k)-\left((1+4\rho)/20k\right)\sum_{s}\frac{1}{k_s}+(d/4k)}\right).
\]
Using the choice of the parameters $k_s$,
$
\sum_s\frac1{k_s}
>
\frac65\rho-\frac9{20}-2\varepsilon,
$ we deduce
\begin{equation}
\sum_{q\in Q}\chi(q)
\ll |Q|^{1-1/2k}\left(p^{(4/20k)\rho-(9/20k)-\left((1+4\rho)/20k\right)((6/5)\rho-9/20-2\varepsilon)+(d/4k)}\right).
\end{equation}
Comparing the estimate with the trivial bound of the character sum, we obtain 
\[
|Q|\geq p^{(4/10)\rho-(9/10)-\left((1+4\rho)/10\right)((6/5)\rho-9/20-2\varepsilon)+(d/2 )}.
\]
Since each $J_i$ has size at least $p^\rho$, the set $Q$ has size at least $p^{(d-1)\rho}$. Substituting this into the above we get a quadratic polynomial
\[
96\rho^2+(200d-292)\rho+(171-100d)>0.
\]
Taking the positive root of the above quadratic equation, we get 
\[
\rho>
\frac{-50d+73 + 5\sqrt{100d^2-196d+49}}{48}=\rho_d''.
\]
\end{proof}
\begin{remark}
$\rho''_3=3/8=.375,\,\rho''_4\approx.417,\,\rho''_5\approx.438,\,\rho''_6\approx.451 \text{ and so on.}$ Moreover, $\rho''_d\rightarrow1/2$ as $d \rightarrow \infty$. Thus, for every fixed d, the above result gives a nontrivial short character sum estimate.
\end{remark}
\begin{remark}
The method does not extend to sub-lattices of higher co-dimension $(>1)$, since in those cases the optimal size of each interval exceeds $p^{1/2}$. Consequently, the same Burgess-type amplification used above no longer obtains a power-saving bound for the character sum.
\end{remark}
\begin{remark}
In general, for any binary homogeneous non-singular (i.e.de-homogenized polynomial has no multiple roots) cubic form, one should also obtain a non-trivial upper bound as in Corollary \ref{MC} whenever the lengths of the intervals are between $p^{3/8+\varepsilon}$ and $p^{1/2}$. 

Suppose first that the cubic form splits completely  into distinct linear factors $(x-\lambda_1y)(x-\lambda_2y)(x-\lambda_3y)$ over $\mathbb{F}_{p}$. In this case, the argument is analogous to that of Case~2 in \cite[Theorem~11]{Cha} for binary quadratic forms. We follow the same strategy by amplifying the intervals $I$ and $J$ using the auxiliary intervals
\[
I_0=\Bigl[1,\frac14 p^{3/8}\Bigr]
\quad\text{and}\quad
K=[1,p^\kappa],
\]
where $\kappa=\varepsilon/4$. We then apply H\"older's inequality to get three factors, as before. The second factor is the key term, which can be estimated using multiplicative energy estimates of $E(R,\,T,\,S)$ where, 
\begin{equation*}
\begin{split}
& R=\{(x-\lambda_1y,x-\lambda_2y,x-\lambda_3y):x\in I,\ y\in J\},\\
& T=\{(x_1-\lambda_1y_1,x_1-\lambda_2y_1, x_1-\lambda_3 y):x_1,y_1\in I_0\},\\
& S=\{(s,s,s):s\in K\}
\end{split}
\end{equation*}
and $ R,\,T,\, S \subset \mathbb{F}_p\times \mathbb{F}_p \times \mathbb{F}_p$. (see \cite[Equations~(6.5),~(6.6) and~(6.7)]{Cha}). The multiplicative energy above can be estimated similarly via lemmas analogous to \cite[Lemma 13,\, 14]{Cha} and following steps same as \cite[Equations~(6.16),~(6.17),~(6.21) and~(6.22)]{Cha}. Finally, choosing $\kappa=3/(2k)$, where $k$ is the exponent arising from H\"older's inequality, yields the desired estimate. 

The remaining case, in which the cubic form factors over $\mathbb{F}_{p}$ into exactly two irreducible factors, can be treated similarly by factoring the form over $\mathbb{F}_{p^2}$ into distinct linear factors and amplifying by the same intervals $I_0$ and $K$ as before. After using H\"older's inequality the key term is estimated by $E(R,T,S)$ as before where,
\begin{equation}
\begin{split}
& R=\{(x-\lambda_1y,x-\omega_2y):x\in I,\ y\in J\},\\
& T=\{(x_1 -\lambda_1 y_1, x_1-\omega_2 y_1): x_1, y_1\in I_0 \},\\
& S=\{(s,s):s\in K\}
\end{split}
\end{equation}
and $\lambda_1\in \mathbb{F}_p,\,\omega_2\in \mathbb{F}_{p^2}$ and $\omega_3\in \mathbb{F}_{p^2}$ are roots of the de-homogenized cubic form. Finally, with the same choice of $\kappa=3/(2k)$ we get a power saving upper bound in character sum for $p^{3/4+\varepsilon}<|I|,|J|<p^{1/2}$.
\end{remark}

{\bf Acknowledgments.} The author thanks Prof. Stephan Baier for the helpful discussions. The author is grateful to Ramakrishna Mission Vivekananda Educational and Research Institute for providing an excellent research environment. This research was supported by a CSIR Ph.D. fellowship under file number 09/0934(13170)/2022-EMR-I.

\end{document}